\documentclass{amsart}

\usepackage{amssymb,mathtools,dsfont,nicefrac}
\usepackage[pdftitle={An explicit rational structure for real semisimple Lie algebras}, pdfauthor={Holger Kammeyer}, pdfsubject={Mathematics}, pdfkeywords={semisimple Lie algebra, rational basis, Chevalley basis}]{hyperref}
\usepackage[lite]{amsrefs}
\usepackage{microtype}
\usepackage{graphicx}
\usepackage{verbatim}
\usepackage[all]{xy}
\usepackage{enumerate}

\newtheorem{Theorem}{Theorem}[section]

\newtheorem{Lemma}[Theorem]{Lemma}
\newtheorem{Proposition}[Theorem]{Proposition}

\theoremstyle{remark}

\newtheorem{Remark}[Theorem]{Remark}
\theoremstyle{definition}
\newtheorem{Definition}[Theorem]{Definition}

\newcommand*{\MRref}[2]{ \href{http://www.ams.org/mathscinet-getitem?mr=#1}{MR \textbf{#1}}}

\newcommand*{\nb}{\nobreakdash}

\newcommand*{\ima}{\textup i}
\newcommand*{\fr}[1]{\mathfrak{#1}}
\newcommand*{\ad}{\textup{ad}}
\newcommand*{\rank}{\textup{rank}}
\newcommand*{\abs}[1]{|#1|}
\newcommand*{\sgn}{\textup{sgn}}

\newcommand*{\Z}{\mathbb Z}
\newcommand*{\R}{\mathbb R}
\newcommand*{\C}{\mathbb C}

\title{An explicit rational structure for real semisimple Lie algebras}
\author{Holger Kammeyer}

\subjclass[2010]{17B20}

\keywords{semisimple Lie algebra, rational basis, Chevalley basis}

\address{Holger Kammeyer\\
Mathematisches Institut\\
Universit\"at Bonn\\
Endenicher Allee 60\\
53115 Bonn\\
Germany\\
kammeyer@math.uni-bonn.de}

\begin{document}


\maketitle

\begin{abstract}
We construct a convenient basis for all real semisimple Lie algebras by means of an adapted Chevalley basis of the complexification.  It determines rational and in fact half-integer structure constants which we express only in terms of the root system and the involution that defines the real structure.
\end{abstract} 

\section{Introduction}

Let \(\fr{g}\) be a complex semisimple Lie algebra with Cartan subalgebra \(\fr{h} \subset \fr{g}\) and Killing form \(B\).  Denote the root system by \(\Phi(\fr{g},\fr{h}) \subset \fr{h}^*\).  Given a root \(\alpha \in \Phi(\fr{g},\fr{h})\), let \(\fr{g}_\alpha \subset \fr{g}\) be its root space and let \(t_\alpha \in \fr{h}\) be the corresponding \emph{root vector} which is defined by \(B(t_\alpha, h) = \alpha(h)\) for all \(h \in \fr{h}\).  Set \(h_\alpha = \frac{2 t_\alpha}{B(t_\alpha, t_\alpha)}\) and for a choice of simple roots \(\Delta(\fr{g},\fr{h}) = \{ \alpha_1, \ldots, \alpha_l\} \subset \Phi(\fr{g},\fr{h})\), set \(h_i = h_{\alpha_i}\).

\begin{Definition}
\label{def:chevalleybasis}
A \emph{Chevalley basis} of \((\fr{g}, \fr{h})\) is a basis \(\mathcal{C} = \{x_\alpha, h_i: \alpha \in \Phi(\fr{g}, \fr{h}), 1 \leq i \leq l \}\) of \(\fr{g}\) with the following properties.
\begin{enumerate}[(i)]
\item \label{item:chevalleybasis:xxh} \(x_\alpha \in \fr{g}_\alpha\) and \([x_\alpha, x_{-\alpha}] = -h_\alpha\) for each \(\alpha \in \Phi(\fr{g}, \fr{h})\).
\item \label{item:chevalleybasis:calphabeta}  For pairs of roots \(\alpha, \beta \in \Phi(\fr{g}, \fr{h})\) such that \(\alpha + \beta \in \Phi(\fr{g}, \fr{h})\), let the constants \(c_{\alpha, \beta} \in \C\) be determined by \([x_\alpha, x_\beta] = c_{\alpha, \beta}x_{\alpha + \beta}\).  Then \(c_{\alpha, \beta} = c_{-\alpha, -\beta}\).
\end{enumerate}
\end{Definition}

This definition appears in \cite{Humphreys:Lie}*{p.\,147}.  In 1955 C.\,Chevalley constructed such a basis and showed that the structure constants are integers \cite{Chevalley:Simple}*{Th\'eor\`eme 1, p.\,24}.

\begin{Theorem}
\label{thm:chevalley}
The structure constants of a Chevalley basis \(\mathcal{C}\) are as follows.
\begin{enumerate}[(i)]
\item \([h_i, h_j] = 0\) for \(i,j = 1, \ldots, l\).
\item \([h_i, x_\alpha] = \langle \alpha, \alpha_i \rangle x_\alpha\) for \(i = 1, \ldots, l\) and \(\alpha \in \Phi(\fr{g},\fr{h})\).
\item \([x_\alpha, x_{-\alpha}] = -h_\alpha\) and \(h_\alpha\) is a \(\Z\)-linear combination of the elements \(h_1, \ldots, h_l\).
\item \label{item:chevalley:calphabeta} \(c_{\alpha, \beta} = \pm (r+1)\) where \(r\) is the largest integer such that \(\beta - r \alpha \in \Phi(\fr{g},\fr{h})\).
\end{enumerate}
\end{Theorem}

As is customary we have used the notation \(\langle \beta, \alpha \rangle = \frac{2 B(t_\beta, t_\alpha)}{B(t_\alpha, t_\alpha)} \in \Z\) with \(\alpha, \beta \in \Phi(\fr{g},\fr{h})\) for the Cartan integers of \(\Phi(\fr{g},\fr{h})\).  The \(\Z\)-span \(\fr{g}(\Z)\) of such a basis is obviously a Lie algebra over \(\Z\) so that tensor products with finite fields can be considered.  Certain groups of automorphisms of these algebras turn out to be simple.  With this method Chevalley constructed infinite series of finite simple groups in a uniform way. For \(\fr{g}\) exceptional he also exhibited some previously unknown ones \cite{Carter:Simple}*{p.\,1}.

But Theorem~\ref{thm:chevalley} states way more than the mere existence of a basis with integer structure constants.  Up to sign, it gives the entire multiplication table of \(\fr{g}\) only in terms of the root system \(\Phi(\fr{g},\fr{h})\).  The main result of this article will be an analogue of Theorem~\ref{thm:chevalley} for any \emph{real} semisimple Lie algebra \(\fr{g}^0\).  To make this more precise, let \(\fr{g}^0 = \fr{k} \oplus \fr{p}\) be a Cartan decomposition of \(\fr{g}^0\) determined by a Cartan involution \(\theta\).  Let \(\fr{h}^0 \subset \fr{g}^0\) be a \(\theta\)-stable Cartan subalgebra such that \(\fr{h}^0 \cap \fr{p}\) is of maximal dimension.  Consider the complexification \((\fr{g}, \fr{h})\) of \((\fr{g}^0,\fr{h}^0)\).  The complex conjugation \(\sigma\) in \(\fr{g}\) with respect to \(\fr{g}^0\) induces an involution of the root system \(\Phi(\fr{g},\fr{h})\).  We will construct a real basis \(\mathcal{B}\) of \(\fr{g}^0\) with (half-)integer structure constants.  More than that, we compute the entire multiplication table of \(\fr{g}^0\) in terms of the root system \(\Phi(\fr{g},\fr{h})\) and its involution induced by \(\sigma\).  For the full statement see Theorem~\ref{thm:structureconstants}.

The idea of the construction is as follows.  We pick a Chevalley basis \(\mathcal{C}\) of \((\fr{g},\fr{h})\) and for \(x_\alpha \in \mathcal{C}\) we consider twice the real and twice the negative imaginary part, \(X_\alpha = x_\alpha + \sigma(x_\alpha)\) and \(Y_\alpha = \ima(x_\alpha - \sigma(x_\alpha))\), as typical candidates of elements in \(\mathcal{B}\).  It is clear that \(\sigma(x_\alpha) = d_\alpha x_{\alpha^\sigma}\) for some \(d_\alpha \in \C\) where \(\alpha^\sigma\) denotes the image of \(\alpha\) under the action of \(\sigma\) on \(\Phi(\fr{g},\fr{h})\).  But to hope for simple formulas expanding \([X_\alpha, X_\beta]\) as linear combination of other elements \(X_\gamma\) we need to adapt the Chevalley basis \(\mathcal{C}\) to get some control on the constants \(d_\alpha\).  A starting point is the following lemma of D.\,Morris \cite{Morris:QForms}*{Lemma 6.4, p.\,480}.  We state it using the notation we have established so far.  Let \(\tau\) be the complex conjugation in \(\fr{g}\) with respect to the compact form \(\fr{u} = \fr{k} \oplus \ima \fr{p}\).

\begin{Lemma}
\label{lemma:morris}
There is a Chevalley basis \(\mathcal{C}\) of \((\fr{g},\fr{h})\) such that for all \(x_\alpha \in \mathcal{C}\)
\begin{enumerate}[(i)]
\item \label{item:morris:tau} \(\tau(x_\alpha) = x_{-\alpha}\),
\item \label{item:morris:sigma} \(\sigma(x_\alpha) \in \{\pm x_{\alpha^\sigma}, \pm \ima x_{\alpha^\sigma} \}\).
\end{enumerate}
\end{Lemma}

In fact Morris proves this for any Cartan subalgebra \(\fr{h} \subset \fr{g}\) which is the complexification of a general \(\theta\)-stable Cartan subalgebra \(\fr{h}^0 \subset \fr{g}^0\).  With our special choice of a so-called \emph{maximally noncompact} \(\theta\)-stable Cartan subalgera \(\fr{h}^0\), we can sharpen this lemma.  We will adapt the Chevalley basis \(\mathcal{C}\) to obtain \(\sigma(x_\alpha) = \pm x_{\alpha^\sigma}\) (Proposition~\ref{prop:existencesigmatau}) and we will actually determine which sign occurs for each root \(\alpha \in \Phi(\fr{g},\fr{h})\) (Proposition~\ref{prop:signofalpha}).  By means of a Chevalley basis of \((\fr{g}, \fr{h})\) thus adapted to \(\sigma\) and \(\tau\) we will then obtain a version of Theorem~\ref{thm:chevalley} over the field of real numbers (Theorem~\ref{thm:structureconstants}).  We remark that a transparent method of consistently assigning signs to the constants \(c_{\alpha, \beta}\) has been proposed by Frenkel-Kac \cite{Frenkel-Kac:Basic}.

A notable feature of the basis \(\mathcal{B}\) is that part of it spans the nilpotent algebra \(\fr{n}\) in an Iwasawa decomposition \(\fr{g}^0 = \fr{k} \oplus \fr{a} \oplus \fr{n}\).  In fact, a variant of \(\mathcal{B}\) is the disjoint union of three sets spanning the Iwasawa decomposition (Theorem~\ref{thm:basisiwasawa}).  For all Iwasawa \(\fr{n}\)-algebras we obtain integer structure constants whose absolute values have upper bound six.  Invoking the classification of complex semisimple Lie algebras, we improve this bound to four (Theorem~\ref{thm:nilpotentupperbound}).  

The outline of sections is as follows.  Section~\ref{sec:restrictedroots} fixes notation and recalls the basic concepts of restricted roots and the Iwasawa decomposition of real semisimple Lie algebras.  Section~\ref{sec:adapted} carries out the adaptation procedure for a Chevalley basis as we have indicated.  In Section~\ref{sec:integral} we prove the Chevalley-type theorem for real semisimple Lie algebras and outline the applications to Iwasawa decompositions.  The material in this article is part of the author's doctoral thesis \cite{Kammeyer:L2-invariants}.

\section{Restricted roots and the Iwasawa decomposition}
\label{sec:restrictedroots}

Again let \(\fr{g}^0\) be a real semisimple Lie algebra with Cartan decomposition \(\fr{g}^0 = \fr{k} \oplus \fr{p}\) determined by a Cartan involution \(\theta\).  There is a maximal abelian \(\theta\)\nb-stable subalgebra \(\fr{h}^0 \subseteq \fr{g}^0\), unique up to conjugation, such that \(\fr{a} = \fr{h}^0 \cap \fr{p}\) is maximal abelian in \(\fr{p}\) \cite{Helgason:Symmetric}*{pp.\,259 and 419--420}.  The dimension of \(\fr{a}\) is called the \emph{real rank} of \(\fr{g}^0\), \(\rank_\R\,\fr{g}^0 = \dim_\R \fr{a}\).  Given a linear functional \(\alpha\) on \(\fr{a}\), let
\[ \fr{g}^0_\alpha = \left\{ x \in \fr{g}^0 \colon [h,x] = \alpha(h)x \text{ for each } h \in \fr{a} \right\}. \] 
If \(\fr{g}^0_\alpha\) is not empty, it is called a \emph{restricted root space} of \((\fr{g}^0, \fr{a})\) and \(\alpha\) is called a \emph{restricted root} of \((\fr{g}^0, \fr{a})\).  Let \(\Phi(\fr{g}^0, \fr{a})\) be the set of restricted roots.  The Killing form \(B^0\) of \(\fr{g}^0\) restricts to a Euclidean inner product on \(\fr{a}\) which carries over to the dual \(\fr{a}^*\).  Within this Euclidean space the set \(\Phi(\fr{g}^0, \fr{a})\) forms an abstract root system.  Note two differences to the complex case.  On the one hand, the root system \(\Phi(\fr{g}^0, \fr{a})\) might not be reduced.  This means that given \(\alpha \in \Phi(\fr{g}^0, \fr{a})\), it may happen that \(2 \alpha \in \Phi(\fr{g}^0, \fr{a})\).  On the other hand, restricted root spaces will typically not be one-dimensional.

Now choose positive roots \(\Phi^+(\fr{g}^0, \fr{a})\).  Then define a nilpotent subalgebra \(\fr{n} = \oplus\, \fr{g}^0_\alpha\) of \(\fr{g}^ 0\) by the direct sum of all restricted root spaces of positive restricted roots.  We want to call it an Iwasawa \(\fr{n}\)-algebra.  The real semisimple Lie algebra \(\fr{g}^0\) is the direct vector space sum of a compact, an abelian and a nilpotent subalgebra, \(\fr{g}^0 = \fr{k} \oplus \fr{a} \oplus \fr{n}\).  The possible choices of positive restricted roots exhaust all possible choices of Iwasawa \(\fr{n}\)-algebras in the decomposition.  Their number is thus given by the order of the Weyl group of \(\Phi(\fr{g}^0, \fr{a})\).  Let \(\fr{g} = \fr{g}^0_\C\) be the complexification.  Then \(\fr{h} = \fr{h}^0_\C\) is a Cartan subalgebra of \(\fr{g}\).  It determines the set of roots \(\Phi(\fr{g},\fr{h}) \subseteq \fr{h}^*\).  Let \(B=B^0_\C\) be the complexified Killing form.  Let \(\fr{h}_\R \subset \fr{h}\) be the real span of the root vectors \(t_\alpha\) for \(\alpha \in \Phi(\fr{g},\fr{h})\).  It is well-known that the restriction of \(B\) turns \(\fr{h}_\R\) into a Euclidean space with orthogonal decomposition \(\fr{h}_\R = \fr{a} \oplus \ima (\fr{h}^0 \cap \fr{k})\).  In what follows, we will need various inclusions as indicated in the diagram
\[
\begin{xy}
\xymatrix@R=15pt@C=25pt{
& \fr{a} \ar[d]^i \ar[dl]_k \ar[r] & \fr{g}^0 \ar[d]^l\\
\fr{h}_\R \ar[r]^j & \fr{h} \ar[r] & \fr{g}.
}
\end{xy}
\]
The compatibility \(l^*B = B^0\) is clear.  Let \(\Sigma = \left\{ \alpha \in \Phi(\fr{g}, \fr{h}) \colon i^* \alpha \neq 0 \right\}\) be the set of roots which do not vanish everywhere on \(\fr{a}\).  The terminology ``restricted roots'' is explained by the following two facts  \cite{Helgason:Symmetric}*{pp.\,263 and 408}.
\begin{enumerate}[(i)]
\item \label{item:restrictedroots1} We have \(\Phi(\fr{g}^0, \fr{a}) = i^* \Sigma\).
\item \label{item:restrictedroots2} For each \(\beta \in \Phi(\fr{g}^0, \fr{a})\), we have \(\fr{g}_\beta^0 = ( \bigoplus\limits_{\substack{i^* \alpha = \beta\\\alpha \in \Sigma}} \fr{g}_\alpha ) \cap \fr{g}^0\).
\end{enumerate}
Statement (\ref{item:restrictedroots1}) says in particular that each \(\alpha \in \Sigma\) takes only real values on \(\fr{a}\).  In fact, \(j^* \Phi(\fr{g}, \fr{h})\) is a root system in \(\fr{h}_\R^*\) and the restriction map \(i^*\) translates to the orthogonal projection \(k^*\) onto \(\fr{a}^*\).

\section{Adapted Chevalley bases}
\label{sec:adapted}

Recall that \(\sigma\) and \(\tau\) denote the complex anti-linear automorphisms of \(\fr{g}\) given by conjugation with respect to \(\fr{g}^0 = \fr{k} \oplus \fr{p}\) and the compact form \(\fr{u} = \fr{k} \oplus \ima \fr{p}\), respectively.  Evidently \(\theta = l^* (\sigma \tau)\) so that \(\sigma \tau\) is the unique complex linear extension of \(\theta\) from \(\fr{g}^0\) to \(\fr{g}\) which we want to denote by \(\theta\) as well. Since \(\sigma\), \(\tau\) and \(\theta\) are involutive, \(\sigma\) and \(\tau\) commute.  Choose positive roots \(\Phi^+(\fr{g}, \fr{h})\) such that \(i^*\Phi^+(\fr{g}, \fr{h}) = \Phi^+(\fr{g}^0, \fr{a}) \cup \{0\}\) and let \(\Delta(\fr{g}, \fr{h}) \subset \Phi^+(\fr{g}, \fr{h})\) be the set of simple roots.  For \(\alpha \in \Phi(\fr{g}, \fr{h})\) let \(h_\alpha = \frac{2}{B(t_\alpha, t_\alpha)} t_\alpha\) and set \(h_i = h_{\alpha_i}\) for the simple roots \(\alpha_i \in \Delta(\fr{g}, \fr{h})\) where \(1 \leq i \leq l = \rank_\C(\fr{g})\).  Let \(\alpha^\sigma, \alpha^\tau, \alpha^\theta \in \fr{h}^*\) be defined by \(\alpha^\sigma (h) = \overline{\alpha(\sigma(h))}\), \(\alpha^\tau(h) = \overline{\alpha(\tau(h))}\) and \(\alpha^\theta (h) = \alpha(\theta(h))\) where \(\alpha \in \fr{h}^*\), \(h \in \fr{h}\).  If \(\alpha \in \Phi(\fr{g}, \fr{h})\) and \(x_\alpha \in \fr{g}_\alpha\), then
\[{[}h,\sigma(x_\alpha)] = \sigma([\sigma(h), x_\alpha]) = \sigma(\alpha(\sigma(h))x_\alpha) = \overline{\alpha(\sigma(h))}\sigma(x_\alpha)\]
so that \(\sigma(x_\alpha) \in \fr{g}_{\alpha^\sigma}\) and similarly for \(\tau\) and \(\theta\).  Thus in this case \(\alpha^\sigma\), \(\alpha^\tau\) and \(\alpha^\theta\) are roots.  Since \(\fr{h}_\R = \fr{a} \oplus \ima (\fr{h}^0 \cap \fr{k})\), we see that \(\alpha^\tau = -\alpha\) for each \(\alpha \in \Phi(\fr{g},\fr{h})\).  We adopt a terminology of A.\,Knapp \cite{Knapp:Beyond}*{p.\,390} and call a root \(\alpha \in \Phi(\fr{g}, \fr{h})\) \emph{real} if it is fixed by \(\sigma\), \emph{imaginary} if it is fixed by \(\theta\) and \emph{complex} in all remaining cases.  Note that \(\alpha^\sigma = -\alpha\) if and only if \(\alpha\) is imaginary.  A real root vanishes on \(\fr{h}^0 \cap \fr{k}\), thus takes only real values on \(\fr{h}_0\).  An imaginary root vanishes on \(\fr{a}\), thus takes purely imaginary values on \(\fr{h}^0\).  A complex root takes mixed complex values on \(\fr{h}^0\).  The imaginary roots form a root system \(\Phi_{\ima \R}\) \cite{Helgason:Symmetric}*{p.\,531}.  The complex roots \(\Phi_\C\) and the real roots \(\Phi_\R\) give a decomposition of the set \(\Sigma = \Phi_\C \cup \Phi_\R\) which restricts to the root system \(i^*\Sigma = \Phi(\fr{g}^0, \fr{a})\).  Let \(\Delta_0 = \Delta(\fr{g}, \fr{h}) \cap \Phi_{\ima \R}\) be the set of simple imaginary roots and let \(\Delta_1 = \Delta(\fr{g}, \fr{h}) \cap \Sigma\) be the set of simple complex or real roots.

Recall Definition~\ref{def:chevalleybasis}, Theorem~\ref{thm:chevalley} and Lemma~\ref{lemma:morris} of the introduction.  Our goal is to prove the following refinement of Lemma~\ref{lemma:morris}.

\begin{Proposition}
\label{prop:existencesigmatau}
There is a Chevalley basis \(\mathcal{C} \!=\! \{x_\alpha, h_i\colon \!\alpha \in \Phi(\fr{g}, \fr{h}), 1 \leq i \leq l \}\) of \((\fr{g}, \fr{h})\) such that
\begin{enumerate}[(i)]
\item \label{item:existencesigmatau:tauadapted} \(\tau(x_\alpha) = x_{\alpha^\tau} = x_{-\alpha}\) for each \(\alpha \in \Phi(\fr{g}, \fr{h})\),
\item \label{item:existencesigmatau:sigmaadapted}
\(\sigma(x_\alpha) = \pm x_{\alpha^\sigma}\) for each \(\alpha \in \Phi(\fr{g}, \fr{h})\) and\newline
\(\sigma(x_\alpha) = + x_{\alpha^\sigma}\) for each \(\alpha \in \Phi_{\ima \R} \cup \Delta_1\).
\end{enumerate}
\end{Proposition}

\begin{Remark}
A.\,Borel \cite{Borel:CliffordKlein}*{Lemma 3.5, p.\,116} has built on early work by F. Gantmacher \cite{Gantmacher:Canonical} to prove a lemma which at least assures that \(\sigma(x_\alpha) = \pm x_{\alpha^\sigma}\) for all \(\alpha \in \Phi(\fr{g},\fr{h})\).  But Borel's method only works for a \emph{maximally compact} \(\theta\)-stable Cartan subalgebra \(\fr{h}^0\) of \(\fr{g}^0\), which is one that has intersection with \(\fr{k}\) of maximal dimension.  We have made the opposite choice of a \emph{maximally noncompact} \(\theta\)-stable Cartan subalgebra \(\fr{h}^0\) that has intersection with \(\fr{p}\) of maximal dimension.
\end{Remark}

We will say that a Chevalley basis \(\mathcal{C}\) is \emph{\(\tau\)-adapted} if it satisfies (\ref{item:existencesigmatau:tauadapted}) and \emph{\(\sigma\)-adapted} if it satisfies (\ref{item:existencesigmatau:sigmaadapted}) of the proposition.  We prepare the proof with the following lemma.

\begin{Lemma}
\label{lemma:satake}
There is a unique involutive permutation \(\omega \colon \Delta_1 \rightarrow \Delta_1\) and there are unique nonnegative integers \(n_{\beta \alpha}\) with \(\alpha \in \Delta_1\) and \(\beta \in \Delta_0\) such that for each \(\alpha \in \Delta_1\)
\begin{enumerate}[(i)]
\item \label{item:satake:alphatheta} \( \alpha^\theta = -\omega(\alpha) - \sum \limits_{\beta \in \Delta_0} n_{\beta \alpha} \beta\)
\item \label{item:satake:nomega} \(n_{\beta\, \omega(\alpha)} = n_{\beta \alpha}\) and
\item \label{item:satake:extension} \(\omega\) extends to a Dynkin diagram automorphism \(\omega \colon \Delta(\fr{g},\fr{h}) \rightarrow \Delta(\fr{g},\fr{h})\).
\end{enumerate}
\end{Lemma}

Part (\ref{item:satake:alphatheta}) is due to I.\,Satake \cite{Satake:Compactifications}*{Lemma~1, p.\,80}.  As an alternative to Satake's original proof, A.\,L.\,Onishchik and E.\,B.\,Vinberg suggest a slightly differing argument as a series of two problems in \cite{Onishchik-Vinberg:Lie}*{p.\,273}.  We will present the solutions because they made us observe the additional symmetry (\ref{item:satake:nomega}) which will be important in the proof of Proposition~\ref{prop:existencesigmatau}.  Part (\ref{item:satake:extension}) can be found in the appendix of \cite{Onishchik:Real}*{Theorem~1, p.\,75}, which was written by J.\,\v{S}ilhan.

\begin{proof}
Let \(C\) be an involutive \((n\!\times\!n)\)-matrix with nonnegative integer entries.  It acts on the first orthant \(X\) of \(\R^n\), the set of all \(v \in \R^n\) with only nonnegative coordinates.  We claim that \(C\) is a permutation matrix.  Since \(C\) is invertible, every column and every row has at least one nonzero entry.  Thus we observe \(\abs{Cv}_1 \geq \abs{v}_1\) for all \(v \in X\).  Suppose the \(i\)-th column of \(C\) has an entry \(c_{ji} \geq 2\) or a second nonzero entry.  Then the standard basis vector \(\varepsilon_i \in X\) is mapped to a vector of \(L^1\)-norm at least \(2\).  But that contradicts \(C\) being involutive.

Let \(\alpha \in \Delta_1\).  Then \(\alpha^\theta\) is a negative root, so we can write
\[ \alpha^\theta = - \textstyle \sum\limits_{\gamma \in \Delta_1} n_{\gamma \alpha} \gamma - \sum\limits_{\beta \in \Delta_0} n_{\beta \alpha} \beta\]
with nonnegative integers \(n_{\gamma \alpha}\) and \(n_{\beta \alpha}\).  Consider the transformation matrix of \(\theta\) acting on \(\fr{h}^*\) with respect to the basis \(\Delta(\fr{g}, \fr{h})\).  In terms of the decomposition \(\Delta(\fr{g}, \fr{h}) = \Delta_1 \cup \Delta_0\) it takes the block form
\[
\left(\begin{array}{c|c}
-n_{\gamma \alpha} & 0 \\ \hline
-n_{\beta \alpha} & \mathbf{1}
\end{array}\right)
\]
with \(\mathbf{1}\) representing the \(\abs{\Delta_0}\)-dimensional unit matrix.  The block matrix squares to a unit matrix.  For the upper left block we conclude that \((n_{\gamma \alpha})\) is a matrix \(C\) as above and thus corresponds to an involutive permutation \(\omega \colon \Delta_1 \rightarrow \Delta_1\).  This proves (\ref{item:satake:alphatheta}).  For the lower left block we conclude that \(n_{\beta \alpha} = \sum_{\delta \in \Delta_1} n_{\beta \delta} n_{\delta \alpha} = n_{\beta \omega(\alpha)}\) because \((n_{\delta \alpha})\) is the aforementioned permutation matrix, so \(n_{\delta \alpha} = 1\) if \(\delta = \omega(\alpha)\) and \(n_{\delta \alpha} = 0\) otherwise.  This proves (\ref{item:satake:nomega}).
\end{proof}

\begin{proof} (of Proposition~\ref{prop:existencesigmatau}.) Pick a Chevalley basis \(\mathcal{C}\) of the pair \((\fr{g}, \fr{h})\).  The proofs of Lemma~\ref{lemma:morris}\,(\ref{item:morris:tau}) by Borel and Morris make reference to the conjugacy theorem of maximal compact subgroups in connected Lie groups.  We have found a more hands-on approach that has the virtue of giving a more complete picture of the proposition:  The adaptation of \(\mathcal{C}\) to \(\tau\) is gained by adjusting the norms of the \(x_\alpha\).  Thereafter the adaptation of \(\mathcal{C}\) to \(\sigma\) is gained by adjusting the complex phases of the \(x_\alpha\).

From Definition~\ref{def:chevalleybasis}\,(\ref{item:chevalleybasis:xxh}) we obtain \(-\frac{2 t_\alpha}{B(t_\alpha, t_\alpha)} = [x_\alpha, x_{-\alpha}] = B(x_\alpha, x_{-\alpha}) t_\alpha\), thus \(B(x_\alpha, x_{-\alpha}) < 0\) because \(B(t_\alpha, t_\alpha) > 0\).  But also \(B(x_\alpha, \tau x_\alpha) < 0\).  Indeed, \((x_\alpha + \tau x_\alpha) \in \fr{k}\) where \(B\) is negative definite, so \(B(x_\alpha + \tau x_\alpha, x_\alpha + \tau x_\alpha) < 0\) and \(B(x_\alpha, x_\beta) = 0\) unless \(\alpha + \beta = 0\).  If constants \(b_\alpha \in \C\) are defined by \(\tau x_\alpha = b_\alpha x_{-\alpha}\) for \(\alpha \in \Phi(\fr{g}, \fr{h})\), we conclude that the \(b_\alpha\) are in fact positive real numbers.  Moreover, \(b_{-\alpha} = b_\alpha^{-1}\) because \(\tau\) is an involution.  We use Definition~\ref{def:chevalleybasis}\,(\ref{item:chevalleybasis:calphabeta}) to deduce \(b_{\alpha + \beta} \!=\! b_\alpha b_\beta\) from \([\tau x_\alpha, \tau x_\beta] \!=\! \tau ([x_\alpha, x_\beta])\) whenever \(\alpha, \beta, \alpha + \beta \in \Phi(\fr{g}, \fr{h})\).  In other words and under identification of \(\alpha\) and \(t_\alpha\), the map \(b\) defined on the root system \(j^*\Phi(\fr{g}, \fr{h})\) extends to a homomorphism from the root lattice \(Q = \Z (j^*\Delta(\fr{g}, \fr{h}))\) to the multiplicative group of positive real numbers.  We replace each \(x_\alpha\) by \(\frac{1}{\sqrt{b_\alpha}} x_\alpha\) and easily check that we obtain a Chevalley basis with unchanged structure constants that establishes (\ref{item:existencesigmatau:tauadapted}).

Now assume that \(\mathcal{C}\) is \(\tau\)-adapted.  It is automatic that \(\sigma(x_\beta) = + x_{\beta^\sigma} = x_{-\beta}\) for each \(\beta \in \Phi_{\ima \R}\) because for each imaginary root \(\beta\) the root space \(\fr{g}_\beta\) lies in \(\fr{k} \otimes \C\) \cite{Helgason:Symmetric}*{Lemma 3.3\,(ii), p.\,260}.  But \(\fr{k} \otimes \C\) is the fixed point algebra of \(\theta\), so the assertion follows from (\ref{item:existencesigmatau:tauadapted}) and \(\sigma = \tau \theta\).  We define constants \(u_\alpha \in \C\) by \(\theta(x_\alpha) = u_\alpha x_{\alpha^\theta}\) for \(\alpha \in \Phi(\fr{g}, \fr{h})\).  As we have just seen, \(u_\alpha = 1\) if \(\alpha\) is imaginary.  In general, the \(\tau\)-adaptation effects \(\sigma(x_\alpha) = \overline{u_\alpha} x_{\alpha^\sigma}\) and \(\overline{u_\alpha} = u_{-\alpha}\) because \(\sigma = \tau \theta = \theta \tau\).  Note
\[-u_\alpha u_{-\alpha} h_{\alpha^\theta} = [u_\alpha x_{\alpha^\theta}, u_{-\alpha} x_{-\alpha^\theta}] = [\theta(x_\alpha), \theta(x_{-\alpha})] = -\theta(h_\alpha) = -h_{\alpha^\theta},\]
so \(u_{-\alpha} = u_\alpha^{-1}\) and \(\abs{u_\alpha} = 1\).  From \(\theta^2(x_\alpha) = x_\alpha\) we get \(u_{\alpha^\theta} = u_\alpha^{-1} = u_{-\alpha} \ (\ast)\).  Next we want to discuss the relation of \(u_\alpha\) and \(u_{\omega(\alpha)}\) for \(\alpha \in \Delta_1\).  First assume that for a given two-element orbit \(\{\alpha, \omega(\alpha)\}\) the integers \(n_{\beta \alpha}\) of Lemma~\ref{lemma:satake} vanish for all \(\beta \in \Delta_0\).  A notable case where this condition is vacuous for all \(\alpha \in \Delta_1\), is that of a \emph{quasi-split} algebra \(\fr{g}^0\) when \(\Delta_0 = \emptyset\).  From \(n_{\beta \alpha} = 0\) we get \(\omega(\alpha)^\theta = -\alpha\).  Thus \(u_{\omega(\alpha)} = u_{-\omega(\alpha)^\theta} = u_\alpha\) by means of (\(\ast\)).  Now assume there is \(\beta_0 \in \Delta_0\) such that \(n_{\beta_0 \alpha} > 0\).  From Lemma~\ref{lemma:satake}\,(\ref{item:satake:alphatheta}) and (\ref{item:satake:nomega}) we get that \(-\omega(\alpha)^\theta = \alpha + \sum_{\beta \in \Delta_0} n_{\beta \alpha} \beta\) is the unique decomposition of \(-\omega(\alpha)^\theta\) as a sum of simple roots.  It is well-known that this sum can be ordered as \(-\omega(\alpha)^\theta = \alpha_1 + \cdots + \alpha_k\) such that all partial sums \(\gamma_i = \alpha_1 + \cdots + \alpha_i\) are roots.  Thus \(x_{-\omega(\alpha)^\theta} = \textstyle \prod_{i=1}^{k-1} \!c_{\alpha_{i+1}, \gamma_i}^{\;-1} \ \ad(x_{\alpha_k}) \cdots \ad(x_{\alpha_2}) (x_{\alpha_1})\).
For one \(i_0\) we have \(\alpha_{i_0} = \alpha\) and the remaining \(\alpha_i\) are imaginary.  Hence by \((\ast)\)
\begin{align*}
u_{\omega(\alpha)} x_{-\omega(\alpha)} &= \theta (x_{-\omega(\alpha)^\theta}) =  \textstyle \prod\limits_{i=1}^{k-1} c_{\alpha_{i+1}, \gamma_i}^{\;-1} \ u_\alpha \: \ad(x_{\alpha_k^\theta}) \cdots \ad(x_{\alpha_2^\theta}) (x_{\alpha_1^\theta}) = \\
&= \textstyle \prod\limits_{i=1}^{k-1} \frac{c_{\alpha_{i+1}^\theta, \gamma_i^\theta}}{c_{\alpha_{i+1}, \gamma_i}} \ u_\alpha \: x_{-\omega(\alpha)} = \pm u_\alpha x_{-\omega(\alpha)}.
\end{align*}
Here we used that \(c_{\alpha, \beta} = \pm c_{\alpha^\theta, \beta^\theta}\) by Theorem~\ref{thm:chevalley}\,(\ref{item:chevalley:calphabeta}) because \(\theta\) induces an automorphism of the root system \(\Phi(\fr{g},\fr{h})\).  It follows that \(u_{\omega(\alpha)} = \pm u_\alpha\) and the sign depends on the structure constants of the Chevalley basis only.  We want to achieve \(u_{\omega(\alpha)} = + u_\alpha\).  So for all two-element orbits \(\{\alpha, \omega(\alpha)\}\) with \(u_\alpha = - u_{\omega(\alpha)}\), replace \(x_{\omega(\alpha)}\) and \(x_{-\omega(\alpha)}\) by their negatives.  This produces a new \(\tau\)-adapted Chevalley basis \(\{ x_\alpha', h_i \colon \alpha \in \Phi(\fr{g}, \fr{h}) \}\) though some structure constants might have changed sign.  Set \(\theta(x_\alpha') = u_\alpha' x_{\alpha^\theta}'\) for all \(\alpha \in \Phi(\fr{g}, \fr{h})\).  We claim that \(u_{\omega(\alpha)}' = u_\alpha'\) for all \(\alpha \in \Delta_1\). The only critical case is that of an \(\alpha \in \Delta_1\) with \(n_{\beta_0 \alpha} > 0\) for some \(\beta_0 \in \Delta_0\).  But in this case we deduce from Lemma~\ref{lemma:satake} that neither \(-\alpha^\theta\) nor \(-\omega(\alpha)^\theta\) is simple, yet only vectors \(x_{\omega(\alpha)}\), \(x_{-\omega(\alpha)}\) corresponding to simple roots \(\omega(\alpha)\) with \(u_{\omega(\alpha)} = -u_\alpha\) have been replaced.  So still \(u_{\omega(\alpha)}' = u_\alpha'\) if we had \(u_{\omega(\alpha)} = u_\alpha\).  If \(u_{\omega(\alpha)} = -u_\alpha\), the replacement is given by \(x_{\alpha^\theta}' = x_{\alpha^\theta}\) and \(x_{\omega(\alpha)^\theta}' = x_{\omega(\alpha)^\theta}\) as well as \(x_\alpha' = x_\alpha\) whereas \(x_{\omega(\alpha)}' = - x_{\omega(\alpha)}\).  Thus,
\begin{align*}
u_\alpha' x_{\alpha^\theta}' &= \theta (x_\alpha') = \theta(x_{\alpha}) = u_\alpha x_{\alpha^\theta} = u_\alpha x_{\alpha^\theta}' \quad \text{and}\\
u_{\omega(\alpha)}' x_{\omega(\alpha)^\theta}' &= \theta(x_{\omega(\alpha)}') = \theta(-x_{\omega(\alpha)}) = -u_{\omega(\alpha)} x_{\omega(\alpha)^\theta} = - u_{\omega(\alpha)} x_{\omega(\alpha)^\theta}'.
\end{align*}
It follows that \(u_{\omega(\alpha)}' = -u_{\omega(\alpha)} = u_\alpha = u_\alpha'\).  Since \(\Delta(\fr{g}, \fr{h})\) is a basis of \(\fr{h}^*\), there exists \(h \in \fr{h}\) such that \(e^{\alpha(h)} = u_\alpha'\) and \((-\ima) \alpha(h) \in (-\pi,\pi]\) for all \(\alpha \in \Delta(\fr{g}, \fr{h})\).  From \(u_\beta' = 1\) for each \(\beta \in \Delta_0\) we get \(h \in \bigcap_{\beta \in \Delta_0} \textup{ker}(\beta)\) and from \(u_{\omega(\alpha)}' = u_\alpha'\) we get \(\alpha(h) = \omega(\alpha)(h)\) for each \(\alpha \in \Delta_1\).  Thus by Lemma~\ref{lemma:satake} we have for each \(\alpha \in \Delta_1\)
\[\alpha^\theta(h) = -\alpha(h).\]
We remark that since \(\alpha(\theta h) = \alpha(-h)\) holds true for all \(\alpha \in \Delta(\fr{g},\fr{h})\), it follows \(\theta(h) = -h\), so \(h \in \ima \fr{a}\).  Let \(x_\alpha'' = e^{-\frac{\alpha(h)}{2}} x_\alpha'\) for each \(\alpha \in \Phi(\fr{g}, \fr{h})\).  Then Definition~\ref{def:chevalleybasis}\,(\ref{item:chevalleybasis:xxh}) and (\ref{item:chevalleybasis:calphabeta}) hold for the new \(x_\alpha''\).  But so does Proposition~\ref{prop:existencesigmatau}\,(\ref{item:existencesigmatau:tauadapted}) because \(\alpha(h)\) is purely imaginary for each \(\alpha \in \Phi(\fr{g}, \fr{h})\) and because \(\tau\) is complex antilinear.  For \(\alpha \in \Delta_1\) we calculate
\[ \theta(x_\alpha'') = e^{-\frac{\alpha(h)}{2}} \theta(x_\alpha') = e^{-\frac{\alpha(h)}{2}} u_\alpha' x_{\alpha^\theta}' = e^{\frac{\alpha(h)}{2}} e^{\frac{\alpha^\theta(h)}{2}} x_{\alpha^\theta}'' = x_{\alpha^\theta}''. \]
From now on we will work with the basis \(\{x_\alpha'', h_i\colon \alpha \in \Phi(\fr{g}, \fr{h})\}\) and drop the double prime.  We have \(\theta(x_\alpha) = x_{\alpha^\theta}\) for each \(\alpha \in \Phi_{\ima \R} \cup \Delta_1\).  It remains to show \(\theta(x_\alpha) = \pm x_{\alpha^\theta}\) for general \(\alpha \in \Phi(\fr{g}, \fr{h})\).  First let \(\alpha \in \Phi(\fr{g}, \fr{h})^+\) be positive and let \(\alpha = \alpha_1 + \cdots + \alpha_k\) be a decomposition as a sum of simple roots such that for \(1 \leq j \leq k\) the partial sums \(\gamma_j = \alpha_1 + \cdots + \alpha_j\) are roots.  Then we have
\[ x_\alpha = \textstyle \prod\limits_{i=1}^{k-1} c_{\alpha_{i+1}, \gamma_i}^{\;-1} \ \ad(x_{\alpha_k}) \cdots \ad(x_{\alpha_2}) (x_{\alpha_1}).\]
Thus
\[ \theta (x_\alpha) = \textstyle \prod\limits_{i=1}^{k-1} c_{\alpha_{i+1}, \gamma_i}^{\;-1} \ \ad(x_{{\alpha_k}^\theta}) \cdots \ad(x_{{\alpha_2}^\theta}) (x_{{\alpha_1}^\theta}) = \textstyle\prod\limits_{i=1}^{k-1} \frac{c_{{\alpha_{i+1}}^\theta, \gamma_i^\theta}}{c_{\alpha_{i+1}, \gamma_i}} \ x_{\alpha^\theta} = \pm x_{\alpha^\theta}.\]
Finally we compute
\[ {[}\theta(x_{-\alpha}), \pm x_{\alpha^\theta}] = [\theta(x_{-\alpha}), \theta(x_\alpha)] = \theta(h_\alpha) = h_{\alpha^\theta} = [x_{-\alpha^\theta}, x_{\alpha^\theta}], \]
hence \(\theta(x_{-\alpha}) = \pm x_{-\alpha^\theta}\).  We conclude \(\theta(x_\alpha) = \pm x_{\alpha^\theta}\) and \(\sigma(x_\alpha) = \pm x_{\alpha^\sigma}\) for all \(\alpha \in \Phi(\fr{g}, \fr{h})\).
\end{proof}

The constructive method of proof also settles two questions that remain.  Which combinations of signs of the \(c_{\alpha, \beta}\) can occur for a \(\sigma\)- and \(\tau\)-adapted Chevalley basis \(\{x_\alpha, h_\alpha\}\)?  And if we set \(\sigma(x_\alpha) = \sgn(\alpha) x_{\alpha^\sigma}\) for \(\alpha \in \Phi(\fr{g}, \fr{h})\), how can we compute \(\sgn(\alpha) \in \{\pm 1\}\)?  We put down the answers in the following two propositions.

\begin{Proposition}
\label{prop:sigmataucriterion}
A set of Chevalley constants \(\{c_{\alpha, \beta}\colon \alpha + \beta \in \Phi(\fr{g},\fr{h})\}\) of \(\fr{g}\) can be realized by a \(\sigma\)- and \(\tau\)-adapted Chevalley basis if and only if for each two-element orbit \(\{\alpha, \omega(\alpha)\}\) of roots in \(\Delta_1\) with \(n_{\beta_0 \alpha} > 0\) for some \(\beta_0 \in \Delta_0\) we have
\[ \textstyle\prod\limits_{i=1}^{k-1} \frac{c_{{\alpha_{i+1}}^\theta, \gamma_i^\theta}}{c_{\alpha_{i+1}, \gamma_i}} = 1\]
where \(-\omega(\alpha)^\theta = \alpha_1 + \cdots + \alpha_k\) with \(\alpha_i \in \Delta(\fr{g},\fr{h})\) and \(\gamma_i = \alpha_1 + \cdots + \alpha_i \in \Phi(\fr{g},\fr{h})\) for all \(i = 1, \ldots, k\).
\end{Proposition}

\begin{proof}
If the condition on the structure constants holds, take a Chevalley basis of \((\fr{g}, \fr{h})\) which realizes them and start the adaptation procedure of the proof of Proposition~\ref{prop:existencesigmatau}.  Thanks to the condition, the replacement \(x_\alpha \mapsto x_\alpha'\) in the course of the proof is the identity map.  The other two adaptations \(x_\alpha \mapsto \frac{1}{\sqrt{b_\alpha}} x_\alpha\) and \(x_\alpha' \mapsto x_\alpha''\) leave the structure constants unaffected.  Conversely, if \(c_{\alpha, \beta}\) are the structure constants of a \(\sigma\)- and \(\tau\)-adapted Chevalley basis, we compute similarly as in the proof of Proposition~\ref{prop:existencesigmatau} that for each such critical \(\alpha \in \Delta_1\) we have
\[ x_{-\omega(\alpha)} = \theta(x_{-\omega(\alpha)^\theta}) = \textstyle\prod\limits_{i=1}^{k-1} \frac{c_{{\alpha_{i+1}}^\theta, \gamma_i^\theta}}{c_{\alpha_{i+1}, \gamma_i}} x_{-\omega(\alpha)}. \qedhere \]
\end{proof}

In particular, for all quasi-split \(\fr{g}^0\) as well as for all \(\fr{g}^0\) with \(\omega = \textup{id}_{\Delta_1}\) all structure constants of any Chevalley basis of \((\fr{g},\fr{h})\) can be realized by a \(\sigma\)- and \(\tau\)-adapted one.  To compute \(\sgn(\alpha)\) first apply \(\sigma\) to the equation \({[}x_\alpha, x_{-\alpha}] = -h_\alpha\) to get \(\sgn(\alpha) \sgn(-\alpha) [x_{\alpha^\sigma}, x_{-\alpha^\sigma}] = -h_{\alpha^\sigma}, \) so \(\sgn(\alpha) = \sgn(-\alpha)\) for all \(\alpha \in \Phi(\fr{g},\fr{h})\).  Moreover, we get the recursive formula
\label{formula:recursion}
\[ \sgn(\alpha + \beta) = \sgn(\alpha) \sgn(\beta) \frac{c_{\alpha^\sigma, \beta^\sigma}}{c_{\alpha, \beta}} \]
for all \(\alpha, \beta \in \Phi(\fr{g},\fr{h})\) such that \(\alpha + \beta \in \Phi(\fr{g},\fr{h})\).  This follows from applying \(\sigma\) to the equation \({[}x_\alpha, x_\beta] = c_{\alpha, \beta} x_{\alpha + \beta}\).  Since \(\sgn(\alpha) = 1\) for \(\alpha \in \Delta(\fr{g},\fr{h})\) the following absolute version of the recursion formula is immediate.

\begin{Proposition}
\label{prop:signofalpha}
Let \(\{x_\alpha, h_\alpha\}\) be a \(\sigma\)- and \(\tau\)-adapted Chevalley basis of \((\fr{g}, \fr{h})\).  If \(\alpha \in \Phi(\fr{g},\fr{h})^+\), let \(\alpha = \alpha_1 + \cdots + \alpha_k\) with \(\alpha_i \in \Delta(\fr{g},\fr{h})\) and \(\gamma_i = \alpha_1 + \cdots + \alpha_i \in \Phi(\fr{g},\fr{h})\) for all \(i = 1, \ldots, k\).  Then
\[\sgn(\alpha) = \textstyle\prod\limits_{i=1}^{k-1} \frac{c_{{\alpha_{i+1}}^\sigma, \gamma_i^\sigma}}{c_{\alpha_{i+1}, \gamma_i}}.\]
\end{Proposition}

It is understood that the empty product equals one.  Also note that \(c_{\alpha^\sigma, \beta^\sigma} = c_{-\alpha^\sigma, -\beta^\sigma} = c_{\alpha^\theta, \beta^\theta}\).  For completeness we still need to comment on how to find a choice of signs for the \(c_{\alpha, \beta}\) in Theorem~\ref{thm:chevalley}\,(\ref{item:chevalleybasis:calphabeta}) as to obtain some set of Chevalley constants to begin with.  This problem has created its own industry.  One algorithm is given in \cite{Samelson:Lie}*{p.\,54}.  A similar method is described in \cite{Carter:Simple}*{p.\,58}, introducing the notion of \emph{extra special pairs} of roots.  A particularly enlightening approach goes back to I.\,B.\,Frenkel and V.\,G.\,Kac in \cite{Frenkel-Kac:Basic}*{p.\,40}.  It starts with the case of \emph{simply-laced} root systems, which are those of one root length only, then tackles the non-simply-laced case.  An exposition is given in \cite{Kac:Infinite}*{Chapters 7.8--7.10, p.\,105} and also in \cite{deGraaf:Lie}*{p.\,189}.  In this picture the product expression appearing in Propositions~\ref{prop:sigmataucriterion} and \ref{prop:signofalpha} can be easily computed.

\section{An explicit rational structure}
\label{sec:integral}

Pick a \(\sigma\)- and \(\tau\)-adapted Chevalley basis \(\mathcal{C}\) of \((\fr{g},\fr{h})\).  Set \(X_\alpha = x_\alpha + \sigma(x_\alpha)\) and \(Y_\alpha = \ima (x_\alpha - \sigma(x_\alpha))\) for \(\alpha \in \Phi(\fr{g}, \fr{h})\).  Let \(H^1_\alpha = h_\alpha + h_{\alpha^\sigma}\) and \(H^0_\alpha = \ima(h_\alpha - h_{\alpha^\sigma})\).  In other words, \(X_\alpha\), \(H^1_\alpha\) are twice the real part and \(Y_\alpha\), \(H^0_\alpha\) are twice the negative imaginary part of \(x_\alpha\), \(h_\alpha\) in the complex vector space \(\fr{g}\) with real structure \(\sigma\).  Let \(Z_\alpha = X_\alpha + Y_\alpha\).  Let \(\Phi_\C^{+\,*}\) be \(\Phi_\C^+\) with one element from each pair \(\{\alpha, \alpha^\sigma\}\) removed and set \(\Phi_\C^* = \Phi_\C^{+\,*} \cup -\Phi_\C^{+\,*}\).  Here, as always, the plus sign indicates intersection with all positive roots.  Pick one element from each two-element orbit \(\{\alpha, \omega(\alpha)\}\) in \(\Delta_1\) and subsume them in a set \(\Delta_1^*\).  Consider the sets
\begin{gather*}
\mathcal{B}_\R = \{ Z_\alpha \colon \alpha \in \Phi_\R \}, \ \ \mathcal{B}_{\ima \R} = \{ X_\alpha, Y_\alpha \colon \alpha \in \Phi_{\ima \R}^+ \}, \ \ \mathcal{B}_\C = \{ X_\alpha, Y_\alpha \colon \alpha \in \Phi_\C^* \}, \\
\mathcal{H}^1 = \{ H^1_\alpha \colon \alpha \in \Delta_1 \setminus \Delta_1^*\}, \ \ \mathcal{H}^0 = \{ H^0_\alpha \colon \alpha \in \Delta_0 \cup \Delta_1^*\}
\end{gather*}
and let \(\mathcal{B}\) be their union.  Note that for \(\alpha \in \Phi_\R\) we have \(Z_\alpha = X_\alpha\) if \(\sgn(\alpha) = 1\) and \(Z_\alpha = Y_\alpha\) if \(\sgn(\alpha) = -1\).  We agree that \(c_{\alpha, \beta} = 0\) if \(\alpha + \beta \notin \Phi(\fr{g}, \fr{h})\) and \(x_{\alpha} = 0\)  thus \(X_\alpha = Y_\alpha = Z_\alpha = 0\) if \(\alpha \notin \Phi(\fr{g},\fr{h})\).  Since \(\langle \beta, \alpha \rangle\) is linear in \(\beta\), we may allow this notation for all root lattice elements \(\beta \in Q = \Z \Phi(\fr{g},\fr{h})\).

\begin{Theorem}
\label{thm:structureconstants}
The set \(\mathcal{B}\) is a basis of \(\fr{g}^0\) and the subsets \(\mathcal{H}^1\) and \(\mathcal{H}^0\) are bases of \(\fr{a}\) and \(\fr{h}^0 \cap \fr{k}\).  The resulting structure constants lie in \(\frac{1}{2} \Z\) and are given as follows.
\begin{enumerate}[(i)]
\item \label{item:structureconstants:hh} Let \(\alpha, \beta \in \Phi(\fr{g},\fr{h})\).  Then \([H^i_\alpha, H^j_\beta] = 0\) for \(i,j \in \{0,1\}\) and
\item \label{item:structureconstants:hx} \( [H^1_\alpha, X_\beta] = \langle \beta + \beta^\sigma, \alpha \rangle X_\beta, \quad [H^1_\alpha, Y_\beta] = \ \ \langle \beta + \beta^\sigma, \alpha \rangle Y_\beta,\) \newline \indent \([H^0_\alpha, X_\beta] = \langle \beta - \beta^\sigma, \alpha \rangle Y_\beta, \quad \; [H^0_\alpha, Y_\beta] = -\langle \beta - \beta^\sigma, \alpha \rangle X_\beta\).
\item \label{item:structureconstants:zz} Let \(\alpha \in \Phi_\R\). Then \newline
\indent \([Z_\alpha, Z_{-\alpha}] = -\sgn(\alpha) 2 H^1_\alpha\) \newline and \(H^1_\alpha\) is a \(\Z\)-linear combination of elements in \(\mathcal{H}^1\).
\item \label{item:structureconstants:xy} Let \(\alpha \in \Phi_{\ima \R}^+\).  Then \newline \indent \([X_\alpha, Y_{\alpha}] = H^0_\alpha\) \newline
and \(H^0_\alpha\) is a \(\Z\)-linear combination of elements \(H^0_\beta\) for \(\beta \in \Delta_0\).
\item \label{item:structureconstants:xx}Let \(\alpha \in \Phi_\C^*\).  Then \newline
\indent \([X_\alpha, X_{-\alpha}] = -H^1_\alpha, \quad [X_\alpha, Y_{-\alpha}] = -H^0_\alpha, \quad [Y_\alpha, Y_{-\alpha}] \,= H^1_\alpha\) \newline where \(H^1_\alpha\) and \(2H^0_\alpha\) are \(\Z\)-linear combinations in \(\mathcal{H}^1\) and \(\mathcal{H}^0\), respectively.
\item \label{item:structureconstants:mainpart} Let \(\alpha, \beta \in \Phi(\fr{g},\fr{h})\) with \(\beta \notin \{-\alpha, -\alpha^\sigma \}\).  Then \newline
\indent \([X_\alpha, X_\beta] = \ \; c_{\alpha, \beta} X_{\alpha + \beta}  \,+  \sgn(\alpha) c_{\alpha^\sigma, \beta} X_{\alpha^\sigma + \beta},\) \newline
\indent \([X_\alpha, Y_\beta] \; = \ \ c_{\alpha, \beta} Y_{\alpha + \beta}  \; + \sgn(\alpha) c_{\alpha^\sigma, \beta} Y_{\alpha^\sigma + \beta},\) \newline
\indent \([Y_\alpha, Y_\beta] \ \, = -c_{\alpha, \beta} X_{\alpha + \beta} + \sgn(\alpha) c_{\alpha^\sigma, \beta} X_{\alpha^\sigma + \beta}\).
\end{enumerate}
\end{Theorem}

In (\ref{item:structureconstants:mainpart}) there is no reason to prefer \(\alpha\) over \(\beta\).  By anticommutativity we have \(\sgn(\alpha)c_{\alpha^\sigma, \beta} X_{\alpha^\sigma + \beta} = \sgn(\beta) c_{\alpha, \beta^\sigma} X_{\alpha + \beta^\sigma}\) and similarly \(\sgn(\alpha)c_{\alpha^\sigma, \beta} Y_{\alpha^\sigma + \beta} = - \sgn(\beta) c_{\alpha, \beta^\sigma} Y_{\alpha + \beta^\sigma}\).  Of course the basis \(2\mathcal{B}\) has integer structure constants.

\begin{proof}
By construction the set \(\mathcal{B}\) consists of linear independent elements and we have \(\abs{\mathcal{B}} = \dim_\C \fr{g} = \dim_\R \fr{g}^0\).  So \(\mathcal{B}\) is a basis.  Moreover, \(\theta(H^j_\alpha) = (-1)^j H^j_\alpha\) for all \(\alpha \in \Phi(\fr{g},\fr{h})\) so that \(\mathcal{H}^1 \subset \fr{a}\) and \(\mathcal{H}^0 \subset \fr{h}^0 \cap \fr{k}\).  Since \(\dim_\R \fr{a} = \abs{\Delta_1} - \abs{\Delta_1^*}\), these subsets generate.  We verify the list of relations.  Part (\ref{item:structureconstants:hh}) is clear.  Part (\ref{item:structureconstants:hx}) is an easy calculation using \(\langle \beta^\sigma, \alpha^\sigma \rangle = \langle \beta, \alpha \rangle\).  Let \(\alpha \in \Phi_\R\).  Then \(Z_\alpha = X_\alpha\) if \(\sgn(\alpha) = 1\) and \(Z_\alpha = Y_\alpha\) if \(\sgn(\alpha) = -1\).  In the two cases we have \([X_\alpha, X_{-\alpha}] = [2x_\alpha, 2x_{-\alpha}] = -4h_\alpha = -2H^1_\alpha\) and \([Y_\alpha, Y_{-\alpha}] = -4[x_\alpha, x_{-\alpha}] = 2 H^1_\alpha\) so we get the first part of (\ref{item:structureconstants:zz}).  We verify that \(H^1_\alpha\) is a \(\Z\)-linear combination within \(\mathcal{H}^1\) for general \(\alpha \in \Phi(\fr{g},\fr{h})\).  Under the Killing form identification of \(\fr{h}\) with \(\fr{h}^*\) the elements \(t_\alpha \in \fr{h}\) correspond to the roots \(\alpha \in \fr{h}^*\).  The elements \(h_\alpha \in \fr{h}\) correspond to the forms \(\frac{2 \alpha}{B(\alpha, \alpha)} \in \fr{h}^*\) which make up a root system as well, namely the dual root system of \(\Phi(\fr{g},\fr{h})\) with simple roots \(\{ h_\beta \colon \beta \in \Delta(\fr{g},\fr{h}) \}\).  We thus have
\[ \textstyle h_\alpha = \sum\limits_{\gamma \in \Delta_1} k_\gamma h_\gamma + \sum\limits_{\beta \in \Delta_0} k_\beta h_\beta \]
with certain integers \(k_\gamma, k_\beta\) which are either all nonnegative or all nonpositive.  Since \(\beta^\sigma = -\beta\) for \(\beta \in \Delta_0\), we have
\[ H^1_\alpha = h_\alpha + h_{\alpha^\sigma} = \textstyle \sum\limits_{\gamma \in \Delta_1} k_\gamma(h_\gamma + h_{\gamma^\sigma}) = \sum\limits_{\gamma \in \Delta_1} k_\gamma H^1_\gamma. \]
From Lemma~\ref{lemma:satake} we see \(\gamma + \gamma^\sigma = \omega(\gamma) + \omega(\gamma)^\sigma\) and \(B(\omega(\gamma),\omega(\gamma)) = B(\gamma, \gamma)\), so
\[ H^1_\gamma = h_\gamma + h_{\gamma^\sigma} = \textstyle \frac{2 t_\gamma}{B(\gamma,\gamma)} + \frac{2 t_{\gamma^\sigma}}{B(\gamma^\sigma, \gamma^\sigma)} = \frac{2 t_{\gamma + \gamma^\sigma}}{B(\gamma,\gamma)} = \frac{2 t_{\omega(\gamma) + \omega(\gamma)^\sigma}}{B(\omega(\gamma),\omega(\gamma))} = h_{\omega(\gamma)} + h_{\omega(\gamma)^\sigma} = H^1_{\omega(\gamma)}\]
and it follows that
\[  H^1_\alpha = \textstyle\sum\limits_{\gamma \in \Delta_1 \setminus \Delta_1^*} ((1-\delta_{\gamma, \omega(\gamma)})k_{\omega(\gamma)} + k_\gamma) H^1_\gamma \]
with Kronecker-\(\delta\).  This proves the second part of (\ref{item:structureconstants:zz}).  Let \(\alpha \in \Phi_{\ima \R}\).  Then
\[ {[}X_\alpha, Y_\alpha] = [x_\alpha + x_{-\alpha}, \ima(x_\alpha - x_{-\alpha})] = 2 \ima h_\alpha = H^0_\alpha. \]
Since the elements \(h_\alpha\) for \(\alpha \in \Phi_{\ima \R}\) form the dual root system of \(\Phi_{\ima \R}\), we see that \(H^0_\alpha\) is a \(\Z\)-linear combination of elements \(H^0_\beta = 2 \ima h_\beta\) with \(\beta \in \Delta_0\).  This proves (\ref{item:structureconstants:xy}).  To prove (\ref{item:structureconstants:xx}) note first that for each \(\alpha \in \Phi(\fr{g},\fr{h})\) the difference \(\alpha - \alpha^\sigma\) is not a root.  Indeed, if it were, then from the recursion formula on p.\,\pageref{formula:recursion} we would get \(\sgn(\alpha - \alpha^\sigma) = \sgn(\alpha)\sgn(-\alpha^\sigma) \frac{c_{\alpha^\sigma,-\alpha}}{c_{\alpha, -\alpha^\sigma}} = -1\) contradicting Proposition~\ref{prop:existencesigmatau}\,(\ref{item:existencesigmatau:sigmaadapted}) because \(\alpha - \alpha^\sigma = \alpha +\alpha^\theta \in \Phi_{\ima \R}\).  With this remark the three equations are immediate.  It remains to show that \(H^0_\alpha\) is a \(\frac{1}{2}\Z\)-linear combination within \(\mathcal{H}^0\).  From the above decomposition of \(h_\alpha\) as a sum of simple dual roots we get
\[H ^0_\alpha = \ima(h_\alpha - h_{\alpha^\sigma}) = \textstyle\sum\limits_{\gamma \in \Delta_1} k_\gamma H^0_\gamma + \sum\limits_{\beta \in \Delta_0} k_\beta H^0_\beta. \]
We still have to take care of \(H^0_\gamma\) for \(\gamma \in \Delta_1 \setminus \Delta_1^*\).  From Lemma~\ref{lemma:satake} we conclude
\[ h_{\gamma^\sigma} = \textstyle\frac{2}{B(\gamma,\gamma)} t_{\gamma^\sigma} = h_{\omega(\gamma)} + \textstyle\sum\limits_{\beta \in \Delta_0} n_{\beta \gamma} \frac{B(\beta, \beta)}{B(\gamma, \gamma)} h_\beta \]
and the numbers \(m_{\beta \gamma} = n_{\beta \gamma} \frac{B(\beta, \beta)}{B(\gamma, \gamma)}\) are integers.  We thus get
\begin{align*}
H^0_\gamma & = \ima(h_\gamma - h_{\gamma^\sigma}) = \ima (h_\gamma - h_{\omega(\gamma)} - \textstyle\sum\limits_{\beta \in \Delta_0} m_{\beta \gamma} h_\beta ) = \\
&\textstyle = -H^0_{\omega(\gamma)} -2 \ima \sum\limits_{\beta \in \Delta_0} m_{\beta \gamma} h_\beta = -H^0_{\omega(\gamma)} -\sum\limits_{\beta \in \Delta_0} m_{\beta \gamma} H^0_\beta.
\end{align*}
If \(\omega(\gamma) \in \Delta_1^*\), this realizes \(H^0_\gamma\) as a \(\Z\)-linear combination in \(\mathcal{H}^0\). If \(\omega(\gamma) = \gamma\), we obtain \(H^0_\gamma = -\frac{1}{2} \sum\limits_{\beta \in \Delta_0} m_{\beta \gamma} H^0_\beta\) and this is the only point where half-integers might enter the picture.  Part (\ref{item:structureconstants:mainpart}) is an easy application of the recursion formula.
\end{proof}

We construct a slight modification of the basis \(\mathcal{B}\).  It is going to be the union of three sets spanning the fixed Iwasawa decomposition \(\fr{g}^0 = \fr{k} \oplus \fr{a} \oplus \fr{n}\).  We start by discussing the Iwasawa \(\fr{n}\)-algebra.  The observation \(\sigma(\fr{g}_\alpha) = \fr{g}_{\alpha^\sigma}\) allows us to state the decomposition in (\ref{item:restrictedroots2}) on p.\,\pageref{item:restrictedroots2} more precisely as
\[ \fr{g}^0_\beta = \bigoplus\limits_{\alpha \in \Phi_\C^*\colon i^*\alpha = \beta} (\fr{g}_\alpha \oplus \fr{g}_{\alpha^\sigma}) \cap \fr{g}^0 \bigoplus\limits_{\alpha \in \Phi_\R\colon i^*\alpha = \beta} \fr{g}_\alpha \cap \fr{g}^0\]
for each \(\beta \in \Phi(\fr{g}^0, \fr{a})\).  It follows that the set
\[ \mathcal{N} = \left\{ X_\alpha, Y_\alpha, Z_\beta\colon \alpha \in \Phi_\C^{* +}, \,\beta \in \Phi_\R^+ \right\} \subset \mathcal{B} \]
is a basis of \(\fr{n}\).  The structure constants are given in Theorem~\ref{thm:structureconstants}\,(\ref{item:structureconstants:mainpart}) so they are still governed by the root system \(\Phi(\fr{g},\fr{h})\).

\begin{Theorem}
\label{thm:nilpotentupperbound}
Every Iwasawa \(\fr{n}\)-algebra has a basis with integer structure constants of absolute value at most four.
\end{Theorem}

\begin{proof}
From Theorem~\ref{thm:structureconstants}\,(\ref{item:structureconstants:mainpart}) we obtain \(2\abs{c_{\alpha, \beta}}\) as an upper bound of the absolute value of structure constants.  Theorem~\ref{thm:chevalley}\,(\ref{item:chevalley:calphabeta}) and the well-known fact that root strings are of length at most four, tell us that \(c_{\alpha, \beta} \in \pm \{1, 2, 3\}\).  The Chevalley constants \(c_{\alpha, \beta} = \pm 3\) can only occur when \(\fr{g}\) contains an ideal of type \(G_2\).  But \(G_2\) has only two real forms, one compact and one split.  A compact form does not contribute to \(\fr{n}\).  For the split form divide all corresponding basis vectors in \(\mathcal{N}\) by two.  Let \(\alpha\) be the short and \(\beta\) be the long simple root.  Then we have just arranged that the equation \([Z_{2\alpha +\beta}, Z_\alpha] = \pm 3 Z_{3\alpha + \beta}\) gives the largest structure constant corresponding to this ideal.  If \(\fr{g}^0\) happens to have an ideal admitting a complex \(G_2\)-structure, then \(\fr{g}\) has two \(G_2\)-ideals swapped by \(\sigma\).  In that case the corresponding two \(G_2\) root systems are perpendicular.  So one of the two summands in every equation of Theorem~\ref{thm:structureconstants} vanishes and the ideal in \(\fr{g}^0\) does not yield structure constants larger than three either.
\end{proof}

Now we consider the maximal compact subalgebra \(\fr{k}\).  For \(\alpha \in \Phi(\fr{g},\fr{h})\) let \(U_\alpha = X_\alpha + \tau X_\alpha = X_\alpha + X_{-\alpha}\) and similarly \(V_\alpha = Y_\alpha + \tau Y_\alpha = Y_\alpha - Y_{-\alpha}\) as well as \(W_\alpha = Z_\alpha + \tau Z_\alpha = U_\alpha + V_\alpha\).  By counting dimensions we verify
\[ \mathcal{K} = \mathcal{H}^0 \cup \left\{ U_\alpha, V_\alpha, X_\beta, Y_\beta, W_\gamma \colon \alpha \in \Phi_\C^{* +}, \beta \in \Phi_{\ima \R}^+, \gamma \in \Phi_\R^+ \right\} \]
is a basis of \(\fr{k}\).  Thus \(\mathcal{K} \cup \mathcal{H}^1 \cup \mathcal{N}\) is a basis of \(\fr{g}^0 = \fr{k} \oplus \fr{a} \oplus \fr{n}\).  The elements \(U_\alpha\), \(V_\alpha\), \(W_\gamma\) are by construction \(\Z\)-linear combinations of elements in \(\mathcal{B}\).  Conversely, the only elements in \(\mathcal{B}\) which do not lie in \(\mathcal{K} \cup \mathcal{H}^1 \cup \mathcal{N}\) are \(X_{-\alpha}, Y_{-\alpha}\) for \(\alpha \in \Phi_\C^{* +}\) and \(Z_{-\beta}\) for \(\beta \in \Phi_\R^+\).  But for those we have \(X_{-\alpha} = U_\alpha - X_\alpha\), \(Y_{-\alpha} = -V_\alpha + Y_\alpha\) and \(Z_{-\alpha} = \sgn(\alpha) (W_\alpha - Z_\alpha)\).  It follows that the change of basis matrices between \(\mathcal{B}\) and \(\mathcal{K} \cup \mathcal{H}^1 \cup \mathcal{N}\) both have integer entries and determinant \(\pm 1\).  Theorem~\ref{thm:structureconstants} thus gives the following conclusion.

\begin{Theorem}
\label{thm:basisiwasawa}
The set \(\mathcal{K} \cup \mathcal{H}^1 \cup \mathcal{N}\) is a basis of \(\fr{g}^0\) spanning the Iwasawa decomposition \(\fr{k} \oplus \fr{a} \oplus \fr{n}\).  The structure constants lie in \(\frac{1}{2} \Z\).
\end{Theorem}

\noindent {\bf Acknowledgements.} The author is indebted to Thomas Schick for many helpful suggestions and to the Deutsche Forschungsgemeinschaft for funding this research.

\end{document}